\documentclass[12pt,reqno]{amsart}
\usepackage[utf8]{inputenc} 
\usepackage[T1]{fontenc}    
\usepackage{lmodern}  

\usepackage{graphics,color}
\usepackage{amssymb,amsmath,amsthm,amscd}
\usepackage{latexsym,verbatim,graphicx,amsfonts}
\usepackage{hyperref}
\usepackage[mathscr]{euscript}
\usepackage{amsmath, amsthm, amssymb}
\usepackage{dsfont}
\usepackage{enumitem}
\usepackage{mathrsfs}
\usepackage{mathtools}
\usepackage{float}

\theoremstyle{plain}
\newtheorem{theorem}{Theorem}[section]

\newtheorem{corollary}[theorem]{Corollary}
\newtheorem{lemma}[theorem]{Lemma}
\newtheorem{proposition}[theorem]{Proposition}
\theoremstyle{definition}

\newtheorem{definition}[theorem]{Definition}

\theoremstyle{remark}
\newtheorem{rem}[theorem]{Remark}

\DeclareMathOperator{\supp}{supp}

\newcommand{\vertiii}[1]{{\left\vert\kern-0.25ex\left\vert\kern-0.25ex\left\vert #1 
    \right\vert\kern-0.25ex\right\vert\kern-0.25ex\right\vert}}

\font\sstext=ecss1000
\font\sssub=ecss1000 at 7pt
\font\sssubsub=ecss1000 at 5pt

\newfam\ssfam
\textfont\ssfam=\sstext
\scriptfont\ssfam=\sssub
\scriptscriptfont\ssfam=\sssubsub



\begin{document}

\title[Complemented subspaces of a $C(K)$-space]{Complemented subspaces of a $C(K)$-space constructed by Candido}
\author[A. Acuaviva]{Antonio Acuaviva}
\address{School of Mathematical Sciences,
Fylde College,
Lancaster University,
LA1 4YF,
United Kingdom} \email{ahacua@gmail.com}

\date{\today}

\begin{abstract}
We classify the complemented subspaces of $C_0(L \times L)$, where $L$ is an ``exotic'' locally compact Hausdorff space recently constructed by Candido under Ostaszewski's $\clubsuit$-principle.
\end{abstract}

\maketitle

\bigskip
\section{Introduction}

The problem of classifying all complemented subspaces of a given Banach space is notoriously challenging. This difficulty is particularly pronounced in the case of $C(K)$-spaces, i.e. spaces of continuous functions on compact Hausdorff topological spaces $K$. To the best of our knowledge, complemented subspaces have only been classified for a limited number of $C(K)$-spaces. \\

For classical Banach spaces, in the separable setting, the list reduces to $c_0$ \cite{pelczynski1960projections} and $C_0(\omega^\omega)$ \cite{benyamini1978extension}. The classification of complemented subspaces of the former can be further extended to arbitrary cardinalities $c_0(\Gamma)$ \cite{granero1998complemented}.

For other non-separable spaces, a classical result of Lindenstrauss \cite{lindenstrauss1967complemented} establishes that every infinite-dimensional complemented subspace of $\ell_\infty \cong C(\beta \mathbb{N})$ is isomorphic to $\ell_\infty$. More recently, Johnson, Kania and Schechtman \cite{johnson2016closed} also added the spaces of the form $\ell^c_\infty(\Gamma)$ to the list of $C(K)$-spaces with known complemented subspaces. \\

In the non-classical setting, a new line of research has produced numerous exotic locally compact spaces whose complemented subspaces can be classified. These constructions are often initially achieved under additional set-theoretical assumptions.

To name a few, Koszmider originally constructed,  under the Continuum Hypothesis or Martin's Axiom, a scattered locally compact Hausdorff space $K$, such that the only infinite-dimensional complemented subspaces of $C_0(K)$ are $c_0$ and $C_0(K)$ itself \cite{koszmider2005decompositions}. These additional set-theoretical assumptions were later removed by Koszmider and Laustsen \cite{koszmider2021banach}. 

In a similar vein, under $\clubsuit$, Koszmider and Zieliński \cite{koszmider2011complementation} constructed a scattered locally compact Hausdorff space $K$ such that every infini\-te-di\-men\-sio\-nal complemented subspace of $C_0(K)$ is isomorphic to either $c_0$, $C_0(\omega^\omega)$, or $C_0(K)$. Candido \cite{candido2021banach} further developed this construction to create a space $L$ where a straightforward classification of operators on $C_0(\alpha \times L)$ for any ordinal $\alpha$ was achievable. This work enabled him to classify the complemented subspaces of $C_0(\omega \times L)$ and $C_0(\omega^\omega \times L)$. Recently, we have extended the classification of operators to a more general setting, leading to the classification of complemented subspaces in other spaces of the form $C_0(\omega \times K)$ \cite{Acuaviva1}.

A further modification by Candido \cite{candido2024few} led to the construction of another scattered, locally compact Hausdorff space $L$ for which the bounded operators on $C_0(L \times L)$ can be explicitly described. We built on this work to classify the complemented subspaces of $C_0(L \times L)$. Following Candido's notation, we let $\Omega_1$ and $\Omega_2$ be the locally compact Hausdorff spaces defined so that $C_0(\Omega_1)$ and $C_0(\Omega_2)$ are isometrically isomorphic to the symmetric and antisymmetric functions on $C_0(L \times L)$, respectively. Our findings are summarized in the following theorem.

\begin{theorem}[$\clubsuit$]\label{th: main}
    Let $L$ be the locally compact Hausdorff space built by Candido \cite{candido2024few}. Then any complemented subspace of $C_0(L \times L)$ is isomorphic to exactly one of the following: $0$, $\mathbb{R}^n$, $c_0$, $C_0(\omega^\omega)$, $C_0(L)^n$, $C_0(\omega \times L)$, $C_0(\omega^\omega \times L)$, $C_0(\Omega_1)$, $C_0(\Omega_2)$ or $C_0(L \times L)$ for some $n \in \mathbb{N}$.

    Consequently, for $j =1, 2$, the complemented subspaces of $C_0(\Omega_j)$ are precisely: $0$, $\mathbb{R}^n$, $c_0$, $C_0(\omega^\omega)$, $C_0(L)^n$, $C_0(\omega \times L)$, $C_0(\omega^\omega \times L)$ and $C(\Omega_j)$.
\end{theorem}

Before we move on to the proof of our main theorem, we recall the remarkable description of operators on $C_0(L \times L)$ as proven by Candido. Specifically, each bounded operator on $C_0(L \times L)$ can be uniquely expressed as a linear combination of six well-understood types of operators, which we introduce below.

We denote by $I$ the identity operator on $C_0(L \times L)$, and define the \emph{transposition operator} $J: C_0(L \times L) \to C_0(L \times L)$ by
\begin{equation*}
    J(f)(x,y) = f(y,x).
\end{equation*}
Recalling that $C_0(L \times L)$ is isometric to the injective tensor product $C_0(L) \widehat{\otimes}_{\varepsilon} C_0(L)$, we can construct operators induced by the tensor product. These are referred to as \emph{matrix operators}. Specifically, given operators $R_1, R_2, R_3, R_4: C_0(L) \to C_0(L)$ with separable range, we define
\begin{equation*}
    \begin{pmatrix} R_1 & R_3 \\ R_2 & R_4 \end{pmatrix} = R_1 \otimes I + J \circ (R_2 \otimes I) + J \circ (I \otimes R_3) + I \otimes R_4.
\end{equation*}

Finally, we have the \emph{diagonal multiplication operators} $M_g$ and $N_h$ for $g, h \in C_0(L)$, defined as follows:
\begin{equation*}
    M_g(f)(x,y) = g(y)f(x,x), \quad N_h(f)(x,y) = h(x)f(y,y).
\end{equation*}

For more details and a deeper discussion of these operators, we refer to the aforementioned work of Candido, where the following characterization of operators on $C_0(L \times L)$ can also be found.

\begin{theorem}[Candido \cite{candido2024few}]\label{th: Cand1}
Assuming Ostaszewski's $\clubsuit$-principle, there exists a non-metrizable scattered locally compact Hausdorff space \( L \) of height \( \omega \) such that, for every operator \( T : C_0(L \times L) \to C_0(L \times L) \), there exist unique scalars \( p, q \in \mathbb{R} \), unique functions \( g, h \in C_0(L) \), unique separable range operators \( R_1, R_2, R_3, R_4 : C_0(L) \to C_0(L) \), and a unique separable range operator \( S : C_0(L \times L) \to C_0(L \times L) \) such that
\begin{equation*}
    T = p I + q J + M_g + N_h + \begin{pmatrix} R_1 & R_3 \\ R_2 & R_4 \end{pmatrix} + S.
\end{equation*}
\end{theorem}

We give a brief layout of the rest of the paper. In Section \ref{sec: Notation-and-technical}, we introduce the basic notation and technical results required for the proof of Theorem \ref{th: main}. The proof is then presented in Section \ref{sec: thr-complemented-subspaces}.

\bigskip
\section{Notation and Technical Results}\label{sec: Notation-and-technical}

We start by introducing some technical notation and recalling some properties of the space constructed by Candido. We adhere to standard notational conventions and adopt the terminology of Candido \cite{candido2024few} and \cite{candido2021banach} whenever appropriate. For two Banach spaces $X$ and $Y$ we write $X \sim Y$ to mean $X$ and $Y$ are isomorphic, while we use $X \cong Y$ to say that they are isometrically isomorphic. \\

From now on $L$ will always denote the scattered locally compact space built by Candido \cite{candido2024few} under Ostaszewski's $\clubsuit$-principle.  We adhere to the convention of identifying each ordinal $\rho < \omega_1$ with $\rho = \{\alpha < \rho \}$ when appropriate. By Candido's construction, there exists a finite-to-one continuous surjection $\varphi: L \to \omega_1$; for each $\rho < \omega_1$ we denote by $L_\rho = \varphi^{-1}[\rho+1] = \{x \in L: \varphi(x) \in \rho+1\}$, so that $\{L_\rho: \rho < \omega_1\}$ is a clopen cover of $L$.
 
For each clopen subset $V$ of a scattered, locally compact Hausdorff space $K$ we isometrically identify $C_0(V)$ with the subspace of functions that vanish outside $V$. Observe that if $V \subseteq K$ is a clopen subset, we have a decomposition
\begin{equation*}
    C_0(K) \cong C_0(K \backslash V) \oplus C_0(V).
\end{equation*}
As an important case of this phenomenon, for each ordinal $ \rho < \omega_1 $, we isometrically identify $ C_0((L \backslash L_\rho)^2) $ and $ C_0(L^2 \backslash (L \backslash L_\rho)^2) $ with complemented subspaces of $C_0(L \times L) $. Given the pivotal role of these subspaces in characterizing complemented subspaces, we refer to them as $ Z_1^\rho = C_0((L \backslash L_\rho)^2) $ and $ Z_2^\rho = C_0(L^2 \backslash (L \backslash L_\rho)^2)$, so that $C_0(L \times L) \cong Z_1^\rho \oplus Z_2^\rho$. Further, we denote by $\pi_1^\rho, \pi_2^\rho: C_0(L \times L) \to C_0(L \times L)$ the projections onto $Z_1^\rho$ and $Z_2^\rho$ respectively. \\

For convenience, given an operator $T$ in the form of Theorem \ref{th: Cand1}, we will denote $T_1 = p I + qJ$ and $T_2 = M_g + N_h + R + S$, where $R = \begin{pmatrix} R_1 & R_3 \\ R_2 & R_4 \end{pmatrix}$, so that $T = T_1 + T_2$. The idea is that for an operator $T$, the behaviour of $T_1$ can be extracted from its action on $Z_1^\rho$, while if we choose $\rho < \omega_1$ big enough, then the behaviour of $T_2$ on $Z_1^\rho$ is almost trivial and it is controlled inside of $Z_2^\rho$. We formalize this in the next result, which was essentially present in Candido's work.

\begin{proposition}\label{prop: 1}
    Given $T: C_0(L \times L) \to C_0(L \times L)$ as in Theorem \ref{th: Cand1}, there exists $\rho < \omega_1$ such that
\begin{enumerate}[label=(\alph*)]
    \item $R[Z_1^\rho] = \{0\}$ and $R[C_0(L \times L)] \subseteq Z_2^\rho$.  \label{it-prop: R}
    \item $S[Z_1^\rho] = \{0\}$ and $S[C_0(L \times L)]  \subseteq Z_2^\rho$. \label{it-prop: S}
    \item $M_g[C_0(L \times L)] \subseteq C_0(L \times L_\rho) \subseteq Z_2^\rho$ and $N_h[C_0(L \times L)] \subseteq C_0(L_\rho \times L) \subseteq Z_2^\rho$. \label{it-prop: MNhg} 
\end{enumerate}
    In particular, it follows that $T_2[C_0(L \times L)] \subseteq Z_2^\rho$.
\end{proposition}
\begin{proof}
    \ref{it-prop: R} is Proposition 5.6 in \cite{candido2024few}, while \ref{it-prop: S} can be found in the proof of Theorem 1.1 in \cite{candido2024few}. \ref{it-prop: MNhg} follows from the form of $M_g$ and $N_h$, coupled with the fact that, since $g,h \in C_0(L)$, there exists $\rho < \omega_1$ such that $\supp g \subseteq L_\rho$ and $\supp h \subseteq L_\rho$.
\end{proof}

When classifying the complemented subspaces of $C_0(L \times L)$, it is helpful to distinguish between subspaces that are, in a sense, \textit{small}. For convenience, we introduce the following definition.

\begin{definition} 
A closed subspace $Y \subseteq C_0(L \times L)$ is called \emph{small} if it is isomorphic to one of the following spaces: 
\begin{equation*}
    0, \quad \mathbb{R}^n, \quad c_0, \quad C_0(\omega^\omega), \quad C_0(L)^n, \quad C_0(\omega \times L), \quad \text{or} \quad C_0(\omega^\omega \times L).
\end{equation*}
We denote by $\mathcal{X}_s$ the collection of all small subspaces of $C_0(L \times L)$. \end{definition}

Finally, we recall some additional properties already proved by Candido.

\begin{proposition}[Candido]\label{prop-candido-smallness}
\begin{enumerate}[label=(\alph*)]
    \item For any $\rho < \omega_1$, $C_0(L \backslash L_\rho) \sim C_0(L)$, $Z_1^\rho \sim C_0(L \times L)$ and $Z_2^\rho \in \mathcal{X}_{s}$. \label{it-prop: Z-small}
    \item If $Y \in \mathcal{X}_{s}$ and $Z$ is a complemented subspace of $Y$, then $Z \in \mathcal{X}_{s}$.\label{it-prop: small-hereditary} 
    \item If $Y \in \mathcal{X}_{s}$ then $C_0(\omega^\omega \times L) \sim C_0(\omega^\omega \times L) \oplus Y$. In particular, $C_0(L \times L) \sim C_0(L \times L) \oplus Y$.\label{it: L-to-L-absorbent}
    \item For any $\rho < \omega_1$, there exists $\rho < \lambda < \omega_1$ such that $C_0(L_\lambda \backslash L_\rho) \sim C_0(\omega^\omega)$. \label{it- prop: omega-omega-equiv}
\end{enumerate}
\end{proposition}
\begin{proof}
\ref{it-prop: Z-small} $C_0(L \backslash L_\rho) \sim C_0(L)$ follows from \cite[Proposition 5.5]{candido2021banach}, $Z_1^\rho \sim C_0(L \times L)$ by the previous part and the injective tensor product identification while $Z_2^\rho \in \mathcal{X}_{s}$ by \cite[Proposition 6.3]{candido2024few}. \ref{it-prop: small-hereditary} Follows from \cite[Remark 5.8]{candido2021banach}, while \ref{it: L-to-L-absorbent} is elementary. Finally, \ref{it- prop: omega-omega-equiv} follows reasoning as in the proof of \cite[Proposition 5.4]{candido2021banach}.
\end{proof}

We would like to explore now how $Z_1^\rho$ and $Z_2^\rho$ transform under $T_1$. To this end, the following definition will be useful.

\begin{definition}
    Given a subset $V \subseteq L \times L$, we say that $V$ is \emph{symmetric} if $(x,y) \in V$ implies $(y,x) \in V$.
\end{definition}

\begin{rem}\label{rem: 1}
    If $V \subseteq L \times L$ is a symmetric clopen set, then $T_1[C_0(V)] \subseteq C_0(V)$. In particular, for any $\rho < \omega_1$ we have $T_1[Z_1^\rho] \subseteq Z_1^\rho$ and $T_1[Z_2^\rho] \subseteq Z_2^\rho$.
\end{rem} 

If $V \subseteq L \times L$ is a symmetric clopen set, we denote by $C_0^s(V)$ the subspace of symmetric functions on $V$, that is $C_0^s(V) = \frac{1}{2}(I + J)[C_0(V)]$. Similarly, we denote by $C_0^{as}(V)$ the subspace of antisymmetric functions on $V$, in other words $C_0^{as}(V) = \frac{1}{2}(I - J)[C_0(V)]$. Recall from the introduction that, according to this definition, we have  $C_0^s(L \times L) \cong C_0(\Omega_1)$ and $C_0^{as}(L \times L) \cong C_0(\Omega_2)$. Naturally, given any symmetric clopen set $V \subseteq L \times L$ we can decompose $C_0(V) = C_0^{s}(V) \oplus C_0^{as}(V)$.

Moreover, $M_g$ and $N_h$ vanish on the antisymmetric functions, so that if we decompose $f \in C_0(V)$ as $f = f^{s} + f^{as}$, where $f^{s} \in C_0^{s}(V)$ and $f^{as} \in C_0^{as}(V)$, we have
    \begin{equation*}
        M_g(f) = M_g(f^{s}) \text{ and } N_h(f) = N_h(f^{s}).
    \end{equation*}

Using this, we can give the first technical ingredient of the proof.

\begin{lemma} \label{lemma-T-small}
    For every operator $T: C_0(L \times L) \to C_0(L \times L)$ there exists $\rho < \omega_1$ such that $T[Z_2^\rho] \subseteq Z_2^\rho$. In particular, if $T$ is a projection then $T[Z_2^\rho]$ is a complemented subspace of $C_0(L \times L)$ and $T[Z_2^\rho] \in \mathcal{X}_{s}$.
\end{lemma}
\begin{proof}
    From Proposition \ref{prop: 1} we can find $\rho < \omega_1$ such that $T_2[Z_2^\rho] \subseteq Z_2^\rho$, while Remark \ref{rem: 1} guarantees $T_1[Z_2^\rho] \subseteq Z_2^\rho$, and thus $T[Z_2^\rho] \subseteq Z_2^\rho$. 

    Suppose now that $T$ is a projection. Since $T[Z_2^\rho] \subseteq Z_2^\rho$, it is easily checked that $T\pi_2^\rho$ is also a projection and $T[Z_2^\rho] = T \pi_2^\rho [C_0(L \times L)]$, so that $T[Z_2^\rho]$ is a complemented subspace. Since $T[Z_2^\rho]$ can also be seen as a complemented subspace of $Z_2^\rho$, Proposition \ref{prop-candido-smallness} \ref{it-prop: Z-small} and \ref{it-prop: small-hereditary} gives $T[Z_2^\rho] \in \mathcal{X}_{s}$.
\end{proof}

Furthermore, since $T_2$ always maps into $Z_2^\rho$, in the case that $T$ is a projection, one can easily deduce the form of $T_1$. Namely, $T_1$ will either be a trivial projection or the projection onto the symmetric or the antisymmetric functions.

\begin{lemma}\label{lemma: main-part-1}
    Let $T: C_0(L \times L) \to C_0(L \times L)$ be a projection. Then either $T_1 = 0$, $T_1 = \frac{1}{2}(I - J)$, $T_1 = \frac{1}{2}(I + J)$ or $T_1 = I$.
\end{lemma}
\begin{proof}
    This is an easy consequence of Proposition \ref{prop: 1} and Remark \ref{rem: 1}. Indeed, these results together with $T$ being a projection give that $T_1$ is a projection if we regard it as an operator from $Z_1^\rho$ into itself. Since $T_1 = pI + qJ$ and $J^2 = I$, an elementary argument yields the desired possible values for $p$ and $q$.
\end{proof}

The following notation will help us handle the subspaces of symmetric and antisymmetric functions effectively.

\begin{definition}
    For $U \subseteq L \times L$ we define the \emph{symmetric complement} of $U$ by $U^{sc} = \{(y,x): (x,y) \in U\}$.
    The set $U$ is called \emph{extremally antisymmetric} if $(x,y) \in U$ implies $(y,x) \not \in U$.
\end{definition}
\begin{rem}\label{rem: 2}
    For any $U \subseteq L \times L$, $V = U \cup U^{sc}$ is symmetric.
    Moreover, if $U$ is extremally antisymmetric then $U \cap U^{sc} = \emptyset$ and $V = U \sqcup U^{sc}$ satisfies $V \cap \Delta(L^2) = \emptyset$, where $\Delta(L^2) = \{(x, x): x \in L\}$ is the diagonal of $L^2$.
\end{rem}

This provides a straightforward representation of the spaces of symmetric and antisymmetric functions on symmetric sets derived from extremally antisymmetric sets $U$.

\begin{lemma}\label{lmm: sym-and-anti}
    Let $U$ be an extremally antisymmetric clopen set. Then
    \begin{equation*}
        C_0(U) \cong C_0^{s}(U \sqcup U^{sc}) \cong C_0^{as}(U \sqcup U^{sc}).
    \end{equation*}
\end{lemma}
\begin{proof}
    The operators $(I + J): C_0(U) \to C_0^{s}(U \sqcup U^{sc})$ and $(I - J): C_0(U) \to C_0^{as}(U \sqcup U^{sc})$ provide the desired isometric isomorphisms, where the range and domain have been restricted as appropriate.
\end{proof}

We are ready to state the last technical ingredient needed for the proof of Theorem \ref{th: main}.

\begin{lemma}\label{lmm: sim-and-anti}

\begin{enumerate}[label=(\alph*)]
    \item For any $\rho < \omega_1$, we have
    \begin{equation*}
        C_0^s(L \times L) \sim C^s_0((L \backslash L_\rho)^2)  \text{ and }  C_0^{as}(L \times L) \sim C^{as}_0((L \backslash L_\rho)^2). 
    \end{equation*}\label{it: simm-and-anti}
    \item We have
    \begin{equation*}
        C_0^s(L \times L) \sim C_0^s(L \times L) \oplus C_0(\omega^\omega \times L) 
    \end{equation*}
    and
    \begin{equation*}
        C_0^{as}(L \times L) \sim C_0^{as}(L \times L) \oplus C_0(\omega^\omega \times L).
    \end{equation*}
    \label{it: abs-sim-and-anti}
\end{enumerate}

\end{lemma}
\begin{proof}
    We proceed with the proof in the symmetric case, with the antisymmetric case being analogous. We start proving \ref{it: simm-and-anti}. By Proposition \ref{prop-candido-smallness} \ref{it- prop: omega-omega-equiv} we can choose $\lambda < \omega_1$ such that $C_0(L_{\lambda} \backslash L_\rho) \sim C_0(\omega^\omega)$. Decompose
    \begin{equation*}
        C_0^{s}(L \times L) \cong C_0^s((L^2 \backslash (L \backslash L_\rho)^2) \oplus C_0^s(U \sqcup U^{sc}) \oplus C_0^s(V) \oplus C_0^s((L \backslash L_\lambda)^2),
    \end{equation*}
    where
    \begin{equation*}
        U = (L_\lambda \backslash L_\rho) \times (L\backslash L_{\lambda+1})
    \end{equation*}
    and
    \begin{equation*}
        V = [(L_{\lambda+1}\backslash L_\rho) \times (L_\lambda \backslash L_\rho)] \cup [ (L_\lambda \backslash L_\rho) \times (L_{\lambda+1}\backslash L_\rho) ],
    \end{equation*}
    see Figure \ref{fig:1} for an intuitive representation.

 \begin{figure}[ht]
    \centering
    \includegraphics[scale=1]{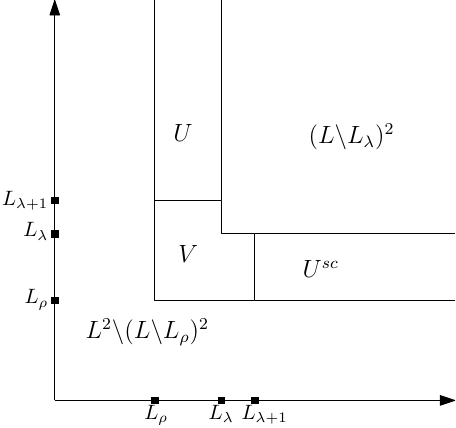}
    \caption{Partition of $L \times L$ in disjoint clopen sets.}
    \label{fig:1}
\end{figure}

    Applying Lemma \ref{lmm: sym-and-anti}, together with our choice of $\lambda$ and Proposition \ref{prop-candido-smallness} \ref{it-prop: Z-small} gives
    \begin{align*}
        C^s_0(U \sqcup U^{sc}) &\cong  C_0(U) \cong C_0(L_\lambda \backslash L_\rho) \hat{\otimes}_\varepsilon C_0(L\backslash L_{\lambda+1}) \\
        &\sim C_0(\omega^\omega ) \hat{\otimes}_\varepsilon C_0(L) \cong C_0(\omega^\omega \times L).
    \end{align*}

On the other hand, $C_0^s((L^2 \backslash (L \backslash L_\rho)^2)$ is a complemented subspace of $C_0((L^2 \backslash (L \backslash L_\rho)^2) = Z_2^\rho$, so that by Proposition \ref{prop-candido-smallness} \ref{it-prop: Z-small} and \ref{it-prop: small-hereditary} we have $C_0^s((L^2 \backslash (L \backslash L_\rho)^2) \in \mathcal{X}_{s}$. Since $ C^s_0(U \sqcup U^{sc}) \sim C_0(\omega^\omega \times L)$, Proposition \ref{prop-candido-smallness} \ref{it- prop: omega-omega-equiv} gives 
\begin{equation*}
     C_0^s((L^2 \backslash (L \backslash L_\rho)^2) \oplus C^s_0(U \sqcup U^{sc}) \sim C^s_0(U \sqcup U^{sc}),
\end{equation*}
and thus
\begin{align*}
     C_0^s(L \times L) \sim  C_0^s(U \sqcup U^{sc}) \oplus C_0^s(V) \oplus C_0^s((L \backslash L_\lambda)^2) \cong C^s_0((L \backslash L_\rho)^2),
\end{align*}
as desired.

To prove \ref{it: abs-sim-and-anti}, we built upon the reasoning before. Observe that $C_0^s(V)$ is a complemented subspace of $C_0(V)$, which in turn is a complemented subspace of $C_0(L^2 \backslash (L \backslash L_\lambda)^2) = Z_2^\lambda$, so that by Proposition \ref{prop-candido-smallness} \ref{it-prop: Z-small} and \ref{it-prop: small-hereditary} we get $C_0^s(V) \in \mathcal{X}_{s}$. Therefore
\begin{equation*}
     C^s_0(U \sqcup U^{sc}) \oplus  C_0^s(V)\sim C^s_0(U \sqcup U^{sc}).
\end{equation*}
It follows that
\begin{equation*}
    C^s_0(L \times L) \sim C_0^s(U \sqcup U^{sc}) \oplus C_0^s((L \backslash L_\lambda)^2) \sim C_0(\omega^\omega \times L)\oplus C^s_0(L \times L),
\end{equation*}
where we have applied part \ref{it: simm-and-anti} to $C_0^s((L \backslash L_\lambda)^2)$.
\end{proof}

Part \ref{it: abs-sim-and-anti} of the previous lemma together with Proposition \ref{prop-candido-smallness} \ref{it- prop: omega-omega-equiv} gives the following.

\begin{corollary}\label{corollary-absorbent}
    Let $Y \in \mathcal{X}_{s}$. Then
    \begin{equation*}
        C^s_0(L \times L) \sim C^s_0(L \times L) \oplus Y \text{ and }  C^{as}_0(L \times L)  \sim C^{as}_0(L \times L) \oplus Y.
    \end{equation*}
\end{corollary}
\bigskip
\section{Proof of Theorem \ref{th: main}.}\label{sec: thr-complemented-subspaces}

We are now ready for the proof of our main result.

\begin{proof}[Proof of Theorem \ref{th: main}]

Each of the spaces is clearly isomorphic to a complemented subspace of $C_0(L \times L)$ and they are pairwise non-iso\-mor\-phic. Indeed, the case $C_0(\omega \times L) \not \sim C_0(\omega^\omega \times L)$ follows from \cite[Corollary 1.2]{candido2021banach} while the rest of the non-trivial cases follow from the discussion in \cite[Section 6]{candido2024few}. Thus, we only need to show that any complemented subspace is isomorphic to one of the claimed ones.

Let $X$ be a complemented subspace of $C_0(L \times L)$ and $T$ be a projection onto $X$. Choose $\rho < \omega_1$ as in Proposition \ref{prop: 1}. From Lemma \ref{lemma: main-part-1}, $T_1$ must take one of four forms, so we divide the proof accordingly. \\

First, assume that $T_1 = 0$ so that $T = T_2 = M_g + N_h + R + S$. Let $f \in C_0(L \times L)$. Since $Tf = T^2f = T(Tf)$ and $Tf \in Z_2^\rho$ by Proposition \ref{prop: 1}, it readily follows that
\begin{equation*}
    X = T[C_0(L \times L)] = T[Z_2^\rho],
\end{equation*}
thus Lemma \ref{lemma-T-small} gives $X \in \mathcal{X}_{s}$. \\

Assume now that $T_1 = \frac{1}{2}(I - J)$ and let $h \in C_0(L \times L)$ be an arbitrary function. Decompose it as $h = f + g$ where $f = \pi_1^\rho h \in Z_1^\rho$ and $g = \pi_2^\rho h \in Z_2^\rho$. Further, write $Tf = f_1 + f_2$ where $f_1 = \pi_1^\rho Tf \in Z_1^\rho$ and $f_2 = \pi_2^\rho Tf \in Z_2^\rho$. 

Observe that since $T_1$ is the projection onto the antisymmetric part, then
\begin{equation*}
    Tf = T_1 f + T_2 f = f^{as} + T_2f,
\end{equation*}
so that $f_1 = \pi_1^\rho Tf = \pi_1^\rho f^{as} \in C_0^{as}((L \backslash L_\rho)^2)$. 
Therefore, $M_g  f_1 = N_h f_1 = 0$, while Proposition \ref{prop: 1} \ref{it-prop: R} and \ref{it-prop: S} guarantee that $R f_1 = S f_1 = 0$. Consequently $T_2 f_1 = 0$ and thus $Tf_1 = T_1 f_1 = f_1$, which gives
\begin{equation*}
    Tf = T^2 f = Tf_1 + T f_2 = f_1 + Tf_2.
\end{equation*}
It follows that
\begin{equation*}
    Th = Tf + Tg = f_1 + T(f_2 + g),
\end{equation*}
so that
\begin{equation*}
    X = T[C_0(L \times L)] \subseteq C_0^{as}((L \backslash L_\rho)^2) + T[Z_2^\rho].
\end{equation*}
The reverse inclusion is easily checked since $T$ acts as the identity on $C_0^{as}((L \backslash L_\rho)^2)$, therefore $X = C_0^{as}((L \backslash L_\rho)^2) + T[Z_2^\rho]$. 
Since $C_0^{as}((L \backslash L_\rho)^2) \subseteq Z_1^\rho$ and $T[Z_2^\rho] \subseteq Z_2^\rho$, we have
\begin{equation*}
    C_0^{as}((L \backslash L_\rho)^2) \cap T[Z_2^\rho] = \{0\}.
\end{equation*}
Moreover, $C_0^{as}((L \backslash L_\rho)^2)$ is a closed subspace, being the image of the projection $\frac{1}{2}(I - J) \pi_1^\rho$, and $T[Z_2^\rho]$ is a closed subspace by Lemma \ref{lemma-T-small}. Thus it follows that
\begin{equation*}
    X = C_0^{as}((L \backslash L_\rho)^2) \oplus T[Z_2^\rho].
\end{equation*}
By Lemma \ref{lmm: sim-and-anti} \ref{it: simm-and-anti} we have
\begin{equation*}
    C_0^{as}((L \backslash L_\rho)^2) \sim C_0^{as}(L \times L),
\end{equation*}
while $T[Z_2^\rho] \in \mathcal{X}_{s}$ by another application of Lemma \ref{lemma-T-small}. Corollary \ref{corollary-absorbent} then gives
\begin{equation*}
    X \sim C_0^{as}(L \times L) \oplus T[Z_2^\rho] \sim C_0^{as}(L \times L) \cong C_0(\Omega_2).
\end{equation*} 

Consider now the case $T_1 = I$, we have that
\begin{equation*}
    X = T[C_0(L \times L)] = T[Z_1^\rho] + T[Z_2^\rho].
\end{equation*}
Since $T_1 = I$ and $\pi_1^\rho T_2 \pi_1^\rho = 0$, it follows that $(T \pi_1^\rho)^2 = T \pi_1^\rho$, so that  $T[Z_1^\rho] = T \pi_1^\rho [C_0(L \times L)]$ is a closed subspace. By Lemma \ref{lemma-T-small}, $T[Z_2^\rho]$ is also a closed subspace. 
Observe that
\begin{equation*}
    T[Z_1^\rho] \cap T[Z_2^\rho] = \{0\}.
\end{equation*}
Indeed, if $f \in T[Z_1^\rho]$, then $f = T\pi^\rho_1 f$, while if $f \in T[Z_2^\rho]$ then $\pi^\rho_1 f = 0$, so that any $f \in T[Z_1^\rho] \cap T[Z_2^\rho]$ satisfies $f = T \pi^\rho_1 f = 0$. It follows that
\begin{equation*}
    X = T[Z_1^\rho] \oplus T[Z_2^\rho].
\end{equation*}
We claim that $T[Z_1^\rho] \sim Z_1^\rho$. Indeed, one can explicitly define an isomorphism between $Z_1^\rho$ and $T[Z_1^\rho]$ by
\begin{equation*}
    A: Z_1^\rho \to T[Z_1^\rho], f \mapsto Tf
\end{equation*}
and
\begin{equation*}
     B: T[Z_1^\rho] \to Z_1^\rho, f \mapsto \pi_1^\rho f.
\end{equation*}
A quick calculation shows $AB  = I_{T[Z_1^\rho]}$ and $BA = I_{Z_1^\rho}$, and thus $T[Z_1^\rho] \sim Z_1^\rho$. Proposition \ref{prop-candido-smallness} \ref{it-prop: Z-small} and \ref{it: L-to-L-absorbent} coupled with Lemma \ref{lemma-T-small} gives
\begin{equation*}
    X \sim Z_1^\rho \oplus T[Z_2^\rho] \sim C_0(L \times L).
\end{equation*}

Lastly, assume that $T_1 = \frac{1}{2}(I + J)$. In this case, $T_2 \pi_1^\rho T_1 \pi_1^\rho = T_2 \pi_1^\rho$ and $\pi_1^\rho T_2 \pi_1^\rho = 0$, so that we still have $(T\pi^\rho_1)^2 = T\pi^\rho_1$. Therefore, arguing as in the previous case, $X = T[Z_1^\rho] \oplus T[Z_2^\rho]$. As before, it is easily checked that in this case $T[Z_1^\rho] \sim C^{s}_0((L \backslash L_\rho)^2)$. Combining Lemmas \ref{lemma-T-small} and \ref{lmm: sim-and-anti} with Corollary \ref{corollary-absorbent} gives $X \sim C^s_0(L \times L)$, finishing the proof.
\end{proof}

\begin{rem} 
Note that it is necessary to consider different cases. When $T_1 = I$ or $T_1 = \frac {1}{2}(I + J)$, the operator $T\pi_1^\rho$ acts as a projection, thus the proof in these two cases is essentially the same. However, this is no longer true when $T_1 = 0$ or $T_1 = \frac{1}{2}(I - J)$, requiring a slightly different argument.
\end{rem}

\bigskip

\noindent\textbf{Acknowledgements.} This paper forms part of the author’s PhD research at Lancaster University under the supervision of Professor N. J. Laustsen. The author extends sincere gratitude to Professor Laustsen for his insightful comments and valuable suggestions on the manuscript’s presentation. The author would also like to extend his gratitude to Professor Leandro Candido for his comments and suggestions on an earlier version of this manuscript and to Bence Horváth for bringing the paper of Candido to his attention.

He also acknowledges with thanks the funding from the EPSRC (grant number EP/W524438/1) that has supported his studies.


\end{document}